\documentclass[10pt,a4paper]{article}
\usepackage[latin1]{inputenc}
\usepackage{amsmath}
\usepackage{graphicx}
\usepackage{amsfonts}
\usepackage{amssymb}
\usepackage{makeidx}
\author{Alex L. Ribeiro-Castro}

\newcommand{\R}{\mathbb R}
\newcommand{\bH}{\mathbb H}
\newcommand{\bI}{\mathbb I}

\newcommand{\sR}{ subRiemannian }
\newcommand{\on}{ orthonormal }

\usepackage{amssymb}
\usepackage{graphicx}
\usepackage{epstopdf}

\newcommand{\g}{\mathfrak{g}}

\newcommand{\fd}{ finite dimensional }

\newcommand{\beq}{\begin{equation}}
\newcommand{\eeq}{\end{equation}}

\newcommand{\C}{\mathbb C}

\newcommand{\Ri}{Riemannian }

\newtheorem{theorem}{Theorem}[section]
\newtheorem{proposition}{Proposition}[section]

\newtheorem{definition}{Definition}[section]

\begin{document}
\begin{center}
\textbf{\begin{large}\textbf{The Chains of Left-invariant CR-structures on SU(2)}\end{large}}\\
Alex Castro and\\
Richard Montgomery,\\
both at the Mathematics Department at UCSC.
\end{center}\medskip

\begin{abstract}We compute the chains associated to the left-invariant CR structures
on the three-sphere.  These structures are characterized by a single
real  modulus $a$.   For the standard structure $a=1$, the chains
are well-known and are   closed curves.   We show that for almost
all other values of the modulus $a$ either two or three types of
chains are simultaneously present : (I) closed curves, (II)
quasi-periodic curves   dense on  two-torii, or (III) chains
homoclinic between closed curves.  For $1 < a < \sqrt{3}$ no curves
of the last type occur.
 A bifurcation occurs at $a = \sqrt{3}$ and from that point on all three types of chains are guaranteed to exist, and exhaust all chains.    The method of proof
  is  to use the Fefferman metric characterization of chains,
  combined with tools from geometric mechanics.  The key to the computation  is a reduced
  Hamiltonian system, similar to Euler's rigid body system, and depending on $a$, which is integrable.
 \end{abstract}

\section{Introduction and Results.}
\label{intro}

The  left-invariant CR structures on the three-sphere $S^3 = SU(2)$ form a
family of CR structures containing the standard structure.
After the standard structure, these form the   most symmetric  CR structures
possible  in  dimension 3.
See Cartan \cite{Cartan}.  The purpose of this note is
to compute  the chains for these   structures.
(Computations of Cartan curvature type invariants for the left-invariant CR structures
can be found in \cite{Cap}.  )

The chains on  a    strictly pseudoconvex CR manifold are a family
of curves on the manifold invariantly associated to its CR
structure.  Chains were defined by Cartan  \cite{Cartan}and further
elucidated by Chern-Moser \cite{ChernMoser}, and  Fefferman
\cite{Fefferman}. Chains  play a role in CR geometry  somewhat
similar to that   of geodesics in \Ri geometry.  The left-invariant
CR structures on $S^3$ are strictly pseudoconvex. Our computation of
the chains for these structures appears here, apparently for the
first time.

The space of left-invariant structures on $S^3 = SU(2)$ modulo
conjugation  is a half-line parameterized by a single  real variable
$a$. Any left-invariant CR structure is conjugate to one of those
presented  in the normal form   below (section \ref{0}, equations (
\ref{normalformA}), (\ref{normalformB}) ).    The standard structure
corresponds to $a =1$. Its  chains  are   obtained by intersecting
$S^3 \subset \C^2$ with complex affine lines in $\C^2$.  (See
\cite{Goldman} for   especially good visual descriptions.)   In
particular all   chains for the standard structure are  closed
curves. Here is our main result:

\begin{theorem} Consider the  left-invariant CR structures on the
three-sphere.  They form a one-parameter space, with parameter $a$
and  $a =1$ corresponding  to the standard structure, as given by  the normal form of
 section \ref{0}, equations
( \ref{normalformA}), (\ref{normalformB}).
Then,  for all but    a discrete set of values of $a$      two types of chains are present:
closed chains  and quasi-periodic chains dense on two-torii.  The curves of each  type are dense   in $S^3$.
A bifurcation occurs at $a = \sqrt{3}$ so that for $a > \sqrt{3}$ a third type of chain occurs,
corresponding to a homoclinic orbit and which accumulates onto a periodic chain ( a geometric circle).
 For  all $a >\sqrt{3}$  all three types of chains: periodic, quasi-periodic, and homoclinic
 are present and every chain is one of these three types.    For $a < \sqrt{3}$ only the closed
 chains and quasi-periodic chains are present.
\end{theorem}

{\bf Remark.} We have left  open the possibility
 that for a finite set of $a \in [1, \sqrt{3}]$ all chains are closed.

 The computations leading to the theorem are based on a construction
 of Fefferman \cite{Fefferman}, refined and generalized by Lee \cite{Lee}
 and Farris \cite{Farris}.   Starting with a strictly
 pseudoconvex CR manifold  $M$ the Fefferman construction yields
 a   circle bundle $S^1 \to X \to M$ together with
 a conformal class of  Lorentzian metrics  on $X$.
 The chains
  are then the projections to $M$ of the
light-like geodesics on $X$.   It follows that  we can look for chains
by solving  Hamiltonian  differential equations.

Once we have the Hamiltonian system for Fefferman's metric, a  simple picture from geometric  mechanics   underlies  this theorem.  For our left-invariant structures  this Hamiltonian  system
is very similar to that of a free rigid body, but with configuration space being
$SU(2) = S^3$ instead of the rotation group $SO(3)$.  Like the rigid body, this Hamiltonian
 system is integrable.  Its  solutions -- the chains -- lie on torii, the  Arnol'd-Liouville torii.
 As in the case of  the rigid body, the non-Abelian symmetry group forces resonances between
 the a priori three  frequencies  on the torii: so that the  torii are
 in fact  two-dimensional, not the expected three dimensions, of $3 = dim(S^3)$.
 When the frequencies are rationally related we get closed  chains.  Otherwise we get the
 quasi-periodic chains.   The phase portrait (figure \ref{fig:2} below) changes
 with $a$ and the bifurcation at $a = \sqrt{3}$  corresponds to the origin turning from an elliptic to
 a hyperbolic fixed point in a bifurcation sometimes known as the Hamiltonian
 figure eight bifurcation.

\subsection{Outline}
There are  five   steps to  the proof of the theorem.  The
  paper is organized along these steps. \\
0.  Find the normal form for the left-invariant structures on $SU(2)$.\\
1. Compute the Fefferman metric on $SU(2) \times S^1$ for the left-invariant CR structures.\\
2. Reduce the  Hamiltonian system
for the Fefferman geodesics by the symmetry group $SU(2) \times S^1$.\\
3.  Integrate the reduced system.\\
4.  Compute the  geometric phases
( holonomies)  relating the full motion to the reduced motion.

We  briefly describe the   methods and ideas involved in each one of
the steps above, and in so doing link  that step to the section in
which it is completed.\\

Step 0. {\bf Finding a normal form. (Section \ref{0} )}
 In section \ref{0} we derive  the normal form (\ref{normalformA}), (\ref{normalformB}) for the  left-invariant CR structures
with single real parameter $a$. This  normal form
 is well-known and standard. Its derivation is   routine. The normal form  can be found
 for example in Hitchin \cite{Hitchin} p. 34,
 and especially the first sentence of the proof of Theorem 10 on p. 99 there.  Hitchin
 provided no derivation of the normal form.
 For completeness  we present the derivation on the normal form in
   section \ref{0}.

Step 1. {\bf Finding the Fefferman metric. (Section \ref{metric})}

In  section \ref{metric} we compute the Fefferman metric associated to our normal forms.
We follow primarily \cite{Lee}.   Inverting this metric yields the Hamiltonian
$H = H_a$ whose solution curves correspond to chains.  \\

Step 2.  {\bf Constructing the reduced dynamics. (Sections
\ref{flows} and \ref{reduced})}The chains for the left-invariant CR
structures are the projections to $S^3$ of the light-like geodesics
for the metrics computed in step 1. These  geodesics are   solutions
to  Hamiltonian systems
  on  $T^*( S^3 \times S^1)$ whose
Hamiltonians  we write    $H = H_a : T^*(S^3 \times S^1) \to \R$.
As with all ``kinetic energy'' Hamiltonians,    $H$ is a   fiber-quadratic function
on the cotangent bundle.     To specify that the geodesics
are light like, we only look at those solutions with    $H=0$.
The  Fefferman metrics are always invariant under   the
circle action.  In our case of left-invariant CR structures
the metrics are also
 invariant under the left action of $S^3 = SU(2)$
(extended in the standard way to the cotangent bundle).
Consequently we can reduce the Fefferman dynamics by the groups
$S^1$ and $SU(2)$.   This reduction is performed in sections \ref{reduced} and  \ref{flows}.  Section \ref{flows} provides
generalities concerning reducing left-invariant flows on Lie groups,
and as such helps to orient the overall discussion.
In section \ref{reduced} we compute the reduced flow.
 In order to perform the reduction
 fix the standard basis
$e_1, e_2, e_3$  for the
$su(2)$.  Write its dual basis, viewed as left-invariant one-forms, as
$\omega^1, \omega^2, \omega^3$.   Write $(g, \gamma)$
for a point of $S^3 \times S^1$ and
$d \gamma$ for the one-form associated to the angular coordinate
$\gamma$.  Any covector $\beta \in T^* _{g, \gamma}(S^3 \times S^1)$
can be expanded as $\beta = M_1  \omega^1 (g)  + M_2  \omega^2  (g)+ M_3  \omega^3 (g)+ P d \gamma$
so  we can write have
$H = H(g, \gamma ;  M_1, M_2, M_3, P)$.  Left-invariance implies
that $H$ does not depend on $g$ or $\gamma$
so we can think of the Hamiltonian as
a function   $ H= H(M_1, M_2, M_3, P)$ on $\R^3 \times \R$.
The Euclidean space     $\R^3 \times \R$  represents $su(2)^* \times \R^{*}$,  the   dual of the Lie algebra of our Lie group, $SU(2) \times S^1$.  Equivalently, $\R^3 \times \R$
is the quotient space
$T^* (S^3 \times S^1)/(S^3 \times S^1)$.   The reduced dynamics is
a flow on this space.
 The coordinate function $P$ is the  momentum map for the
 action of the circle factor and as such is constant along solutions for
 the reduced dynamics.
The function $H$ generates the reduced dynamics:   $\dot M_i = \{ M_i, H \}$
and $\dot P = \{ P, H \} (= 0)$ where $\{ \cdot, \cdot \}$ is the `Lie-Poisson bracket'.
See section \ref{flows}.\\

Step 3.  {\bf Solving the reduced dynamics. (Section
\ref{phase_portrait})} The   phase portrait found in   figures
\ref{fig:1}, and \ref{fig:2} summarizes the reduced dynamics. .  The
computations proceed as follows. The functions  $P$ and $K = M_1 ^2
+ M_2 ^2 + M_3^2$ are Casimirs for the Lie-Poisson structure,
meaning that $\{K, h \} = \{P, h \} = 0$  for any Hamiltonian $h$
used to   generate the reduced dynamics.   The solutions to the
reduced dynamical equations   thus lie on  the  curves formed by the
intersections of the three surfaces $P = const.$,  $K = const.$ and
$H = 0$ in $\R^4 = \R^3 \times \R$. For typical values of these
constants , these curves are closed curves. At special values the
curves may be isolated points, or may be singular, like in the case
of the homoclinic eight (figure \ref{fig:2}).

When $P = 0$ we can solve for the dynamics explicitly.  The corresponding chains
are the left translates of a particular one-parameter subgroup in $G = S^3$.
The case $P \ne 0$ can be reduced to $P =1$ by the following scaling argument.
We have
$H(\lambda M_1,  \lambda M_2, \lambda M_3,  \lambda P) = \lambda^2 H(M_1, M_2, M_3, P)$.
Up on $S^3 \times S^1$ this scaling represents leaving positions alone and scaling
momenta, and hence velocities.  Thus the  reduced solution curves with initial conditions
$(\lambda M_1,  \lambda M_2, \lambda M_3,  \lambda P)$ and those with
initial conditions $(M_1, M_2, M_3, P)$ represent the same geodesics, and so the same
chains, just parameterized differently.  Choosing $\lambda = 1/P$
we can always scale the case $P \ne 0$  to the case $P =1$.  Now we have a single Hamiltonian
$h = H(M_1, M_2, M_3, 1)$ on the standard rigid body phase space $\R^3$.
We represent the surface $h = 0$   as a graph $M_3 = q(M_1, M_2;a)$
 over the $M_1 M_2$ plane, where $q$ is an even quartic function of $M_1, M_2$.
 We form the  solution  curves by intersecting this graph
 with the level sets of $K$.  To simplify the analysis
 we project the resulting curves onto the $M_1 M_2$
 plane.  A critical point analysis of $K$ restricted to the graph  locates the bifurcation value $a = \sqrt{3}$
for the reduced phase portrait as described in theorem 1.

Step 4.  {\bf Geometric phases. (Section \ref{phase})}   We follow
the   idea presented in  the paper \cite{rmontC}   in order to
reconstruct the
  chains in $S^3$ from the reduced solution curves.
  Some mild modifications  are needed to that idea,
  since our   initial   group is $SU(2) \times S^1$
  rather than the group $SO(3)$ of that paper.       Fix $P =1$ and a value of
  $K$ so that the reduced curve $C$ of step 1 is closed.
   The  left action of $SU(2) \times S^1$ on
$T^*( S^3 \times S^1)$ has a momentum map with values in
$su(2)^* \times \R^*$ and solutions (chains) must lie on constant level sets of
this momentum map.    One factor of this momentum map is $P$
from steps 2 and 3 which we have set to $1$.
Upon projecting the   level set onto $T^*S^3$ via the projection $T^* S^3 \times T^*S^1 \to T^* S^3$
  we obtain an embedded $S^3 \subset T^*S^3$ (the graph of a right-invariant one-form)
together with a projection onto the reduced phase space $\R^3 \times \{1 \}$ of step 3.
The inverse image of $C$ under this projection is  a two-torus,
and all the chains whose reduced dynamics is represented by $C$
and whose momentum map has the given fixed value lie on this
two-torus. One angle of
this torus represents  the reduced curve.
 The relevant question  is: as we go once around the reduced curve,
how much does the other angle change?  Call this amount
$\Delta \theta$.   If the value of $\Delta \theta$  is an  irrational multiple of $2 \pi$ then the
chain is not closed and forms one of the quasi-periodic chains of theorem 1, dense on
its   two-torus.
If its value of $\Delta \theta$  is a rational multiple $p/q 2 \pi$
of $2 \pi$   then the chain is closed, corresponding to some
$p,q$ winding on its torus.   With  certain modifications,
the basic integral  formula for $\Delta \theta$ from \cite{rmontC}
is valid.     One term in this formula
 corresponds to a holonomy of a connection, and is termed the
 ``geometric phase'', explaining  the subtitle we have given to this step 4.
 The values of $\Delta \theta$
 depends only on the values of $a$ and $K$
and its dependence is analytic in these variables.  Thus the
proof of the theorem will be complete once we have shown
there is a value of $a$ for which $K \mapsto \Delta \theta (K, a)$
is not constant.

In order to prove non-constancy of  $ \Delta \theta (K, a)$ ,
take $a > \sqrt{3}$ so that the reduced dynamics
has a homoclinic eight.  Denote the value of $K$ on the eight by  $k(a)$.
We show that as   $K \to k(a)$
we have that $\Delta \theta (a, K) \to \infty$.

Steps 0--4 now  completed,   theorem 1 is proved.

\vskip .2cm

{\bf Appendices.} We finish the paper with  two appendices.  In
appendix \ref{app1} we verify  that when $a =1$ the Fefferman
geodesics for the Hamiltonian computed here (eq. \ref{Hamn})
correspond to the well-known chains for the  standard three-sphere.
In appendix 2 we show that the left-invariant CR structures for $a
\ne 1$ correspond to the family of non-embeddable CR structures on
$S^3$ discovered by Rossi, and frequently found in the CR
literature.

{\bf An Open problem.}  We end appendix \ref{app2} with an open
problem inspired by  the Rossi embedding
of $S^3/(\hbox{antipodal map})$ and a conversation with Dan Burns.

\section{A normal form for  the  left-invariant CR structures (step 0). }
\label{0}

\subsection{Preliminaries.  Basic Definitions.}
A contact structure in dimension 3 is   defined
by the vanishing of a one-form $\theta$  having the property
that $\theta \wedge d \theta \ne 0$.   Let $M$ be the underlying 3-manifold
and $TM$ its tangent bundle.    The contact structure is the  field of 2-planes
$\xi= \{ (m, v)  \in TM : \theta(m) (v) = 0 \} \subset TM$.
It is a rank   2 sub-bundle  of the tangent bundle.
The one-form $\theta$ and $f \theta$, for $f \ne 0$ a function, define the same contact structure.

\begin{definition}\label{def1}
A  strictly pseudoconvex CR structure  on a 3-manifold $M$
consists of    contact
structure  $\xi$ on $M$  together with an almost complex structure $J$ defined
  on the contact planes $\xi$.
  \end{definition}

  We will primarily be using the following alternative,
  equivalent definition
  \begin{definition} \label{def2}
A  strictly pseudoconvex CR structure  on a 3-manifold $M$
consists of an oriented contact structure $\xi$ on $M$
together with a conformal equivalence class of metrics
defined on   contact planes $\xi$.
  \end{definition}

 To pass from the first definition to the second,
 we construct the conformal structure from the almost complex structure $J$
 in the standard way.  Namely, the conformal structure
 is determined by knowing what  an orthogonal frame is,
 and we declare   $e, J (e)$ to be such a frame,
  for any nonzero vector $e \in \xi$.
  An alternative to this construction is to
    choose a contact form $\theta$  for the contact structure
    and then construct its associated
  {\it Levi form}
  \beq
  L_{\theta} (v, w) = d \theta (v, J w)
  \label{Leviform}
  \eeq
  which is a quadratic symmetric form on the contact planes.
  The contact condition implies that the Levi form is either negative definite or positive definite.
  If it is negative definite,  replace $\theta$ with $- \theta$ to  make it positive definite.
  We henceforth insist  that $\theta, J$ are taken so the Levi form is positive definite.
  This assumption on $(\theta, J)$ is equivalent to assuming that the orientation on the contact
  planes induced by   $\theta$ and $J$ agree.
  (Note that a choice of
  contact one-form orients the contact planes. )
  The conformal structure associated to $(\theta, J)$
from definition \ref{def1} is generated by the Levi form.
If we change   $\theta \to f \theta$ with $f> 0$
then the Levi form  changes by  $L_{\theta} \to f L_{\theta}$,
showing that this definition of conformal structure is independent of
(oriented) contact form $\theta$.

 To go from definition \ref{def2} to definition \ref{def1},  take any oriented
  orthogonal basis vectors $E_1, E_2$ having the
  same length relative to some metric in the conformal class.
  Define $J$ by  $J(E_1) = E_2, J(E_2) = -E_1$.
    Thus in dimension 3 we can define a CR structure   by
  a contact form $\theta$, defined up to positive scale factor,
  together with  an inner product   on the contact planes $\omega = 0$
  to represent the conformal structure,
  also only defined up to a positive scaling.     Choosing  the scale factor
  of either the contact form or the quadratic form
  fixes the scalar factor of the other  one through  the Levi-form relation, eq. (\ref{Leviform}).

\subsection{The left-invariant case.}
We take $M = S^3$  which  we identify   with
  the Lie group   $SU(2)$ in the standard way, via the action of $SU(2)$
  on $S^3 \subset \C^2$.    A left-invariant CR
structure on $S^3$  is then  given by Lie algebraic data on $su(2)$.
This data consists of a ray in    $su(2)^*$
  representing the left-invariant  contact form  $\theta$ up to positive scale
  and
  a quadratic form on $su(2)$ defined modulo $\theta$, and positive definite
  when restricted to $ker(\theta)$.
       Conjugation on $SU(2)$
maps left invariant CR structures to left-invariant CR structures,
and induces the co-adjoint action on $su(2)^*$.
This action is equivalent, as a representation, to the standard
action of the rotation group $SO(3)$ on $\R^3$
via the  2:1 homomorphism $SU(2) \to SO(3)$.
Consequently, we can rotate the contact form $\theta$
to anti-align with the basis element $\omega_3$.
Thus we take $\theta = - \omega^3$.
The contact planes are then framed  by the left-invariant
vector fields $e_1, e_2 \in su(2)$.    The choice
of $-\omega^3$ is made so that $e_1, e_2$ is the correct orientation of
the plane, as follows from the structure equation
$$d  \omega^3 =  - \omega^1 \wedge \omega^2.$$
This structure equation also proves that the plane    field $-\omega^3  = 0$ is indeed contact,
so that the corresponding CR  structure (no matter   the choice
of $J$)  will be strictly pseudoconvex.
A  quadratic form on the contact plane is   given by
a positive definite  quadratic expression
in $\omega^1, \omega^2$ , that is:
$A ( \omega^1 )^2 + 2B \omega^1 \omega^2 + C (\omega^2)^2$, viewed mod $\omega^3$.      The isotropy group of $\omega^3$
acts by rotations of the contact plane (the $e_1, e_2$ plane).
 A quadratic form can     be diagonalized  by  rotations,
 so upon conjugation by some element of  the   isotropy  subgroup of $\omega^3$
 we can put  the quadratic  form in the diagonal form  $A (\omega^1) ^2 + B (\omega^2) ^2$ with $A, B  > 0$.
 The form is only well-defined up to scale, and we can scale it
 so that $A = 1/B$, i.e the conformal structure
 is that of   $(1/a) (\omega^1) ^2 + a   (\omega^2) ^2$, $a > 0$.
 We have proved the bulk of :
 \begin{proposition} [Normal form]
 \label{normalformprop}
  Every left-invariant CR structure on $S^3$ is conjugate to
 one  whose contact form  is given
 by   \beq
\theta =  - \omega^3     \label{normalformA}
 \eeq
 and whose associated conformal structure is
\beq
L_{\theta} = \frac{1}{a} (\omega^1) ^2 + a  (\omega^2) ^2
\label{normalformB}
\eeq
 The associated almost complex structure $J= J_a$
is defined by
$J(e_1) = {1 \over a} e_2$,  $J(e_2) = - a e_2$.
The structure defined by $a$ is isomorphic to the
structure defined by $1/a$.  As the notation indicates,
the quadratic form $L_{\theta}$ is indeed the Levi-form
associated to $\theta, J$ as per eq. (\ref{Leviform}).
\label{normalform}
\end{proposition}

To see that $J$ in the proposition is correct,
note that the choice $\theta = -\omega^3$ as contact form
induces the orientation $\{ e_1, e_2 \}$ to the contact
planes, and that $\left\lbrace e_1, {1 \over a} e_2 \right\rbrace $ are orthogonal vectors
having the same squared length  ($1/a$)relative to the given metric $L_{\theta}$.
To see that the structure defined by $a$ is isomorphic to the structure defined by
${1 \over a}$ observe that  rotation by 90 degrees converts $(1/a) (\omega ^1) ^2 + a ( \omega^2) ^2$
 to $ a (\omega ^1) ^2 + (1/a) ( \omega^2) ^2$.
 Finally, compute from $d \theta = \omega^1 \wedge \omega^2$
 and the form of $J$ that indeed, the Levi form is the given quadratic
 form $L_{\theta}$.

 \section{Fefferman's metric (step 1).  }
 \label{metric}

When the strictly convex CR structure on $M$
is induced by an embedding $M \subset \C^2$,
Fefferman \cite{Fefferman} constructed a circle bundle $Z \to M$
together with a conformal Lorentzian metric on $Z$ invariantly  associated to
the CR structure.  Farris \cite{Farris}
and then Lee \cite{Lee}  generalized Fefferman's construction
to the case of an abstract strictly pseudoconvex CR structure, i.e. one not necessarily
induced by an embedding into $\C^2$.  In this section we construct the Fefferman metric
for the family of left-invariant CR structures from step 1 (proposition \ref{normalformprop} there).  We most closely follow Lee's presentation.

 We begin with a general construction. Let   $\pi: Z \to M$ be any circle bundle over $M$.
 Fix a contact form $\theta$.  Recall that the Reeb vector field
 associated to $\theta$ is the vector field on $M$ uniquely
 defined by the two conditions
 $$\theta (R) = 1$$
 $$i_R d \theta = 0.$$
 Changing $\theta$ to $g \theta$, $g$ a function, changes
 $R$ to ${1 \over g} R + X_g$ where $X_g$ lies in the contact plane  field
 and is determined pointwise by  a linear equation involving $dg$
 and $d \theta$ which is reminiscent of the equation relating a Hamiltonian
 to its Hamiltonian vector field.  We extend the Levi form (\ref{Leviform})  to all of $TM$
 by insisting that $L_{\theta} (R, v) = 0$ for all $v \in TM$ and continue to
 write   $L_{\theta}$ for this extended  form.
 Let $\sigma$ be any one-form on  $Z$ with the property
 that $\sigma$ is nonzero on the vertical vectors (the kernel of $d \pi$).
 Then
 \beq
 g_{\theta} = \pi^* L_{\theta} + 4  (\pi^* \theta) \odot \sigma
 \label{Feff1}
 \eeq
 is a Lorentzian metric on $Z$. Here $\odot$ denotes the symmetric
 product of one-forms:
 $\theta \odot \sigma = \frac{1}{2}( \theta \otimes \sigma +  \sigma \otimes \theta)$.
\begin{quote}
 {\it The trick needed is a way of  defining  $\sigma$ in terms
 of the contact form, and   $J$,  in such a way that a ``conformal
 change''  $\theta \mapsto g \theta$ of the contact structure
  induces a conformal change of the metric $g_{\theta}$.}
\end{quote}
{\bf Warning.}  Farris and Lee, use a different definition of the symmetric product
$\odot$: their $\theta \odot \sigma$ is twice ours, so that
in their formula for the metric our $4$ is replaced by a $2$.
We have chosen our definition so that, using it,
$(dx + dy)^2 = dx^2 + 2 (dx \odot dy) + dy^2$, where $\theta^2 = \theta \otimes \theta$.

 \subsection{Forming  the circle bundle from the   canonical bundle.   (2,0) forms.}
 \label{circle bundle}

  The circle bundle $Z \to M$ will be a bundle of complex-valued 2-forms,
 defined up to  real scale factor. A choice of contact form $\theta$ on $M$
 induces various   one-forms on $Z$ in a canonical way.
 One of these one-forms will be   the form $\sigma$ needed for the
 Fefferman  metric,  eq. (\ref{Feff1}).
 Here are the main steps leading to the construction of $Z$ and its one-form $\sigma$.

 The complexified contact plane $\xi_{\C} = \xi \otimes \C$
 splits  under $J$ into the holomorphic and anti-holomorphic directions, these being the
 $+i$ and $-i$ eigenspaces of $J$, where $J$ is
  extended from $\xi$ to $\xi_{\C}$ by complex
 linearity.  In the case of 3-dimensional CR manifold, if we start with
 any  non-zero vector field  $E$ tangent to $\xi$,   then
 $Z = E - i JE$ spans the holomorphic direction, while $\bar Z = E + i JE$ spans the anti-holomorphic
 direction.    In our case
 \beq
 Z_a = e_1 - { i \over  a} e_2
 \label{holomorphic}
 \eeq
 is holomorphic, while
 \beq
 \bar Z_a = e_1 +  { i \over  a} e_2
  \label{antiholomorphic}
 \eeq
 is the anti-holomorphic vector field.

 {\bf Remark. Third definition of a 3-dimensional CR manifold.}
Eq. (\ref{holomorphic}) corresponds to yet a third definition of a CR
manifold.
\begin{definition} (CR structure, 3rd time 'round).
 A CR structure on $M^3$
 is a complex line field, i.e. a rank 1 subbundle of the    complexified tangent bundle $TM \otimes {\C}$
 which is nowhere real.
 \end{definition}
  Such a complex  line field is locally spanned by a ``holomorphic'' vector field $Z$ as in eq. (\ref{holomorphic}).
 Writing    $Z = E_1 -  i E_2$ with $E_1, E_2$ real vector fields, we define
 the 2-plane field $\xi$ to be the real span of $E_1, E_2$, and
 we set  $J(E_1) = E_2$, $J(E_2) = - E_1$.
 The  ``strictly pseudoconvex'' condition, which is the condition that
 $\xi$ be contact,  is that  $E_1, E_2$
 together with the Lie bracket
 $[E_1, E_2]$ span the real  tangent bundle $TM$.

  \vskip .2cm

  The almost complex structure $J$
 on  the contact planes of a    CR manifold induces  a splitting of the space of
 complex-valued differential  forms into
 types $\Omega^{p,q}$ similar to  the splitting of forms on a   complex manifolds.
  We declare that a complex valued k-form $\beta$ is of type $(k,0)$ (that is to say ``holomorphic'') if
 $i_{\bar Z} \beta = 0$ for all anti-holomorphic vector fields $\bar Z$.
 In dimension $3$, one only needs to check this equality for a single nonzero such vector field,
 such as $\bar Z$ of eq. ( \ref{antiholomorphic}).

{\bf Our case.} The space of (1,0) forms for the left-invariant structure for the parameter value $a$ is spanned by,
 \beq
 \theta = -\omega^3 \; ; \omega_a  = \omega^1 + i a \omega^2  \, \, \hskip 1.5cm :(1,0) \hbox{ forms }.
 \label{10}
 \eeq
 The (2,0) forms are spanned (over $\C$) by
  \beq
 \theta \wedge \omega_a   \, \, \, \hskip 1.5cm :  (2,0) \hbox{ forms }.
 \label{10}
 \eeq
 In dimension 3  the space of all $(2,0)$ forms, considered pointwise,
 forms a complex line bundle, denoted by $K$ and called the   canonical bundle
  as in  complex differential geometry.   $Z$
 is defined to be the
 ``ray projectivization'' of $K$:
 $$Z =  K \setminus \{ \hbox{ zero section } \} / \R^+.$$

    We next recall  from Lee \cite{Lee}  how a choice of  contact form $\theta$
    determines  the one-form $\sigma$ on $Z$.

 1.  Volume normalization equation.
 Fix the  contact form $\theta$ on $M$.  The volume normalization equation is
 \beq
 \sqrt{-1}\textrm{ }\theta \wedge i_R \zeta \wedge i_R \bar \zeta = \theta \wedge d \theta.
 \label{volnorm}
 \eeq
 The right hand side is the standard volume form defined by a choice of contact
 structure.  On the left-hand side, $R = R_{\theta}$ is the Reeb vector field for
 $\theta$.  The  2-form $\zeta \in \Gamma (K)$, a section of the canonical bundle  is to   viewed as the unknown.
 The equation is quadratic in the unknown since
 multiplying $\zeta$
 by the complex  function $f$ multiplies the left hand side of the volume normalization
 equation by $|f|^2$.  It follows by this scaling
  that there is a   solution, $\zeta_0$ to  the volume normalization
 which is unique up to unit complex multiple $\zeta \mapsto  e^{ i \gamma} \zeta$.

  Said slightly  differently, eq. (\ref{volnorm}) defines a section 
  $$s = s_{\theta}: Z \to K $$
 of the ray bundle $K \to Z$,
 since once we fix the complex phase of $\zeta$, the equation uniquely determines the real scaling factor.
   Fix a solution, which is to say,  a smoothly varying pointwise choice  of solutions
$$\zeta_0: M \to K$$ to
eq. (\ref{volnorm}).  Such a solution choice
 defines  a  global trivialization of $Z$, since we can express
any  point  $z$ of $Z$ (uniquely) as
$$s_{\theta} (z)  = e^{ i \gamma} \zeta_0 (\pi(z)) $$
where $m = \pi(z) \in M$.  Thus the choice $\zeta_0$
induces a   global trivialization:
$$Z \cong M \times S^1.$$
(A more pictorial, equivalent  description of this trivialization of $Z$ is as follows. Form the  ray generated by $\zeta_0 (m)$, which is a point in the circle  fiber $Z_m$, over $m$.
Rotate this ray  by the angle $\gamma$
until you hit the ray $z \in Z_m$, thus associating to $z$ a point $(m, \gamma) \in M \times S^1$).\\
We henceforth use this identification $Z = M \times S^1$
and define a global one-form on $Z$ by
\beq
\zeta (m, \gamma) = e^{ i \gamma} \zeta_0  (m).
\label{zeta}
\eeq

We check now that  the two-form $\zeta$ depends only
on the choice of contact form $\theta$, and so, up to this choice,
is   intrinsic to    $Z$.
  The total space  $K$ of the canonical bundle , like any  total space constructed as a   bundles of $k$-forms, has
  {\it on it} a  canonical $k$-form $\Xi$.  To describe $\Xi$    write a typical point of $K$ as
 $(m, \beta) \in K$, $m \in M$,  $\beta \in \Lambda^{(2,0)}T_x M$.
 Then we can set   $\Xi (x, \beta) = \pi_x ^* \beta$ where
 $\pi: K \to M$ denotes the projection.
 This canonical form, like all such canonical  forms, enjoys the
 {\it  reproducing
 property}  that if $\beta: M \to K$
 is any section, then $\beta ^* \Xi = \beta$.
Let $s= s_{\theta}: Z \to K$ to pull back $\Xi$:
$$ \zeta := s_{\theta} ^* \Xi, \textrm{ a (2,0)-form on $Z$}. $$
The reproducing property shows that, under the global trivialization of $Z$
induced by $\zeta_0$, we have  that $\zeta$ is given by formula (\ref{zeta}) below.\\

{\bf Our case.} Return to the left-invariant situation:   Choosing $\theta = - \omega^3$ we get
 $\theta \wedge d \theta =  - \omega^1 \omega^2 \omega ^3$.  The associated Reeb field is
 \beq
 R = - e_3.
 \label{Reeb}
 \eeq
  Writing $\zeta_0 = g \theta \wedge \omega_a $
 we compute that  $i_R \zeta_0 = g \omega_a$.  Using $\omega_a \wedge \bar \omega_a = -2 i a \omega^1 \wedge \omega^2$ we compute that the left-hand side of the
 volume normalization equation (\ref{volnorm})
 expands out to $-2a |g|^2 \omega^1 \omega^2 \omega ^3$.
 The volume normalization equation (\ref{volnorm})  then   implies that $|g|^2 = 1/2a$.
 Thus
 \beq
 \zeta_a  = \frac{1}{\sqrt{2a} }\theta \wedge \omega_a
 \label{zeta0}
 \eeq
  is a global normalized
 section of $K$.   It
 induces a  global  trivialization of $Z$, as just described, so that we can think of
 $Z$ as  $S^3 \times S^1$.
 With $(m, e^{i \gamma})$ being the ray through
 the (2,0) form $e^{i \gamma} \zeta_a (m)$.
The two-form $\zeta$ on $Z$ is given, under this identification,
 by  this same algebraic relation:
 \beq
 \zeta = e^{i \gamma}   \frac{1}{\sqrt{2a} } \theta \wedge \omega_a
  \label{zeta}
 \eeq
 where we are  not using different symbols to differentiate
 between a form $\beta$ on $M$ and its  pull-backs $\pi^* \beta$ to $Z$.

\begin{proposition}   [Lee: \cite{Lee}, p. 417]

Fix the   contact form $\theta$ for the CR manifold  $M$.
Let $\zeta$ be the induced  one-forms on $Z$
as just described.    Let $R$ be  the Reeb vector field  for $\theta$.

A.  There is a complex valued one-form $\eta$ on $Z$, uniquely determined by
the conditions: .
 \beq
\zeta =\theta\wedge \eta
\label{LeeA}
\eeq
\beq
i_{v} \eta = 0 \hbox{ whenever } \pi_*v = R
\label{LeeB}
\eeq

B.    With $\eta$ as in A, there is a unique real-valued one form $\sigma$ on $Z$
determined by the equations
 \beq
d\zeta = 3i \sigma\wedge \zeta
\label{Lee1}
\eeq
\beq
\label{Lee2}
\sigma\wedge d\eta\wedge\bar{\eta}=Tr(d\sigma)i\sigma\wedge\theta\wedge\eta\wedge\bar{\eta}.
\eeq
The meaning of  $`Tr'=$Trace in this last equation is as follows.
Any solution $\sigma$ to (\ref{Lee1})   has the property that  $d \sigma$  is basic, i.e.
 is the pull-back of a two-form on $M$, which by abuse of notation we  also denote by $d \sigma$.
Any two-form on $M$
can be expressed as   $f d \theta + \theta \wedge \beta$.
Set  $Tr(f d \theta + \theta \wedge \beta) = f$.

C.  The form $\sigma = \sigma(\theta)$ determined by the equations
(\ref{LeeA}, \ref{LeeB}, \ref{Lee1}, \ref{Lee2}) is the form
$\sigma$ appearing in the Fefferman metric $g_{\theta}$ of eq.
(\ref{Feff1}).  If $\theta \mapsto f \theta$, $f > 0$ then the
Reeb  extended Levi form $L_{\theta}$   and $\sigma$ transform in
such a way that
 $g_{f \theta} = f g_{\theta}$, i.e. the conformal class of the  Fefferman metric is indeed invariantly attached
 to the CR structure.

  \end{proposition}

{\sc Remark.}  An equivalent definition of the trace used in eq (\ref{Lee2})  is as follows.
Take a two-form such as $d \sigma$ on $M$,  restrict it   to
 the contact plane and then  use     the Levi form $L_{\theta}$
to raise its  indices and thus  define its trace, $Tr(d \sigma)$.

{\bf  The forms on $Z$ in the  left-invariant case.}  In our left-invariant situation
the   forms $\theta, \zeta$ of the   theorem have    been described above
in equations (\ref{normalformA}), (\ref{zeta}).
They are
$\theta = - \omega^3$,  $\zeta =  \theta \wedge \eta$
with
\beq
\eta =  \frac{1}{\sqrt{2a}} (  e ^{ i \gamma}  \omega_a)
\label{eta}
\eeq
and
$$\omega_a =    (\omega^1 + i a \omega^2) $$
This   $\eta$ is indeed the $\eta$ of part A of the theorem,  since if $V$ is
any vector field on $Z$ satisfying  $\pi_* V = R$
then $i_V \pi^* \eta = i_R \eta = 0$. (Recall we use $\eta$
for $\pi^* \eta$ as forms on $Z$.)

Now we move to the computations of part B of the Proposition for  the one-form $\sigma$.
We compute:
\beq
\sigma =   {{d \gamma} \over {3}} + f \theta \hskip .5 cm , \hskip .5cm  f = \frac{1}{8}( a + 1/a).
\label{sigma}
\eeq
 Here are key  steps along the way of  the computation:
\begin{eqnarray}
d \eta = i d \gamma \wedge \eta + \frac{1}{\sqrt{2a}} e^{i \gamma} d \omega_a \\
= i d \gamma \wedge \eta + \frac{1}{\sqrt{2a}} e^{i \gamma} \theta \wedge (- \omega^2 + i a \omega^1).
\end{eqnarray}
Then
$$ d \zeta = i d \gamma \wedge \zeta $$
It then follows from the first equation in part B of the theorem, and the reality of $\sigma$ that
$$\sigma = \frac{d \gamma}{3} + f \theta$$
for some real function $f$.  We have $Tr(d \sigma) =f$.
Setting $dvol = d \gamma \wedge \theta \wedge \omega^1 \wedge \omega^2$
we compute the right  hand side of  eq. (\ref{Lee2})  to be
$(f/3) dvol$, while its left hand side is equal to
$ [(1/3) (1 + a^2)/2a - f]dvol$.  Setting the two 4-forms equal
and solving for $f$ yields $f = (1 + a^2)/8a$ as claimed.

Returning now to the form of the Fefferman metric, eq. (\ref{Feff1}),
and using $\theta = - \omega^3$
we see that the metric  is given (up to conformality)  by
\beq
ds^2 = \{ {1 \over a} (\omega ^1) ^2 + a (\omega^2) ^2  \}
+  4  \omega^3 \odot
( \frac{1}{8}(a + \frac{1}{a}) \omega^3  - {{d \gamma} \over {3}}  ).
\label{Feff2}
\eeq
Written in terms of the basis  $ \{ e_1, e_2, e_3,  {\partial \over {\partial \gamma}} \}$
this metric is
\beq
g(a)= \left( \begin{array}{cccc}
\frac{1}{a} & 0 & 0 & 0 \\
0 & a & 0 & 0\\
0 & 0 &\frac{1}{2}(a + \frac{1}{a}) & -\frac{2}{3}\\
0 & 0 &- \frac{2}{3} & 0 \end{array} \right).
\label{matrix}
\eeq

\section{Reduced light ray equations (step 2.)}
\label{reduced}
The geodesics for any metric $ds^2 = \Sigma g_{ij} dx^i dx ^j$,
 \Ri or Lorentzian, can be characterized as the  solutions to
 Hamilton's equations for the Hamiltonian defined by
 inverting the metric, and viewing the result as a fiber quadratic
 function on the cotangent bundle:
 \beq
 H(x,p) = {1 \over 2} \Sigma g^{ij} (x) p_i p_j.
 \label{metricH}
 \eeq
 (See for example, \cite{Foundations}, \cite{Arnold}, or \cite{montA}.)
 Here $g^{ij} (x)$ is the matrix pointwise inverse to the matrix
 with entries $g_{ij} (x)$.

 If we are  only interested in    light-like geodesics,   then we restrict
 to   solutions for which  $H = 0$.  It is important that these geodesics are conformally invariant.
 If   $\tilde ds^2 = f ds^2$ is a metric conformal to the original,
 then the corresponding Hamiltonians are related by  $\tilde H = H/ f$ and the two Hamiltionian
 vector fields,  are related  {\it on their common zero level set} $\{ H =  0\}$ by
 $X_{\tilde H} = (1/f) X_H$.  This proportionality of vector fields says
 that the set of light rays for any two conformally related metrics $ds^2, \tilde ds^2$
 are the same  as sets of unparameterized curves.

The Hamiltonian for the Fefferman metric lives on $T^*Z$.
Any covector $p \in T^* _z Z$ can be expanded in the basis
$\omega_1, \omega_2, \omega_3, d \gamma$ dual to the basis
in which the matrix (\ref{matrix}) was computed:
$$p = M_1 \omega_1 + M_2 \omega_2 + M_3 \omega_3 + P d \gamma$$
  The inverse matrix  to  (\ref{matrix})  is
 \begin{equation}
g^{-1}(a)= \left( \begin{array}{cccc}
a& 0 & 0 & 0 \\
0 & {1 \over a} & 0 & 0\\
0 & 0 &0 & -3/2 \\
0 & 0 &-3/2 & - \frac{9}{8}(a + \frac{1}{a})\end{array} \right).
\label{matrixinv}
\end{equation}
It follows that the Fefferman  Hamiltonian for our left-invariant CR structure with parameter $a$ is given by
\begin{equation}
H_a(g, \gamma; M_1,M_2,M_3,P)=\frac{1}{2}\{a M_1 ^2 + {1 \over a}M_2 ^{2} - 3  M_3 P
 - \frac{9}{8}(a + \frac{1}{a})  P^{2} \}.
 \label{Hamn}
\end{equation}

\section{Left-invariant geodesic flows.}
\label{flows}
Our Hamiltonian (\ref{Hamn}, \ref{matrix}) generates the geodesic flow for
a left-invariant (Lorentzian) metric on the Lie group $G = SU(2) \times S^1$.
In this section we review  some  general facts regarding left-invariant geodesic flows,
and   specify to our situation. We refer the reader to \cite{Foundations}, especially chapter 4,  or
\cite{Arnold}, especially Appendix 2,  for  background and more details regarding the material of this
section and the next.

\subsection{Generalities}
Let $Q$ be a manifold.  Let $ds^2$ be  a metric on $Q$ as above. The
geodesic flow for $ds^2$ is encoded  by a   Hamiltonian vector field
$X$   on $T^*Q$ which is defined in terms of the Hamiltonian above
in eq. (\ref{metricH}).
 The vector
field $X$ can be defined by   the
canonical Poisson brackets
$\{  ,  \}$  on $T^*Q$ according to $X[f] = \{ f, H \}$, for
  $f$ any   smooth  function on $T^* Q$.
 It is worth noting that the momentum scaling property  $H (q, \lambda p) = \lambda^2 H(q, p)$,
 for $p \in T^* _q Q$  corresponds to the
fact that the geodesic $\tilde \gamma (t)$ with initial conditions $(q, \lambda p)$ is simply the
same geodesic  $\gamma(t)$ as   represented by the initial conditions $(q, p)$ but just
parameterized  at a different speed: $\tilde \gamma (t) = \gamma (\lambda t)$

Now suppose that   $Q = G$ is  a \fd Lie group
and      the metric is left-invariant, i.e. left translation by any    element of $G$
acts by isometries relative $ds^2$. The left action of $G$ on itself
canonically lifts to $T^*G$, and
left-invariance of the metric  implies that the Hamiltonian  $H$ is left-invariant
under this lifted action.   Write $\g$ for the
Lie algebra of $G$, and $\g^*$ for the dual vector space to
$\g$, which we identify with $T^* _e G$, where $e \in G$ is the identity.
Using the codifferential of left-translation,
we    left-trivialize $T^*G = G \times \g^*$, and use   corresponding
notation $(g, M) \in G \times \g^*$ for points in the trivialized cotangent bundle.   Then the left-invariance of $H$  means that, relative to  this
trivialization we have
$$H(g, M) = H(M)$$
depending  on $M$ alone.

Let $e_a$ be a basis for $\g$, the Lie algebra of $G$, and $\omega^a$ the
corresponding dual basis for $\g^*$.   Then we can expand
$$ M = \Sigma M_a \omega^a$$
and
$$H = {1 \over 2} \Sigma  g^{ab} M_a M_b$$
where $g^{ab}$ is the matrix inverse to the inner product matrix
$g_{ab} = ds^2 (e_a, e_b)$.     We find that
$$\{ M_a, M_b \} = - \Sigma c^d _{ab} M_d$$
where $c^d _{ab}$ are the structure constants of $\g$ relative to the basis
$e_a$.

It follows that the geodesic flow can be pushed down to the quotient space
$(T^* G)/ G = \g^*$, and as such it is represented in coordinates
by
$$ \dot M_a = - \Sigma _ { k, b, r}  g^{ r b} c^k _{a b} M_r M _k $$
We will call these the ``reduced equations", or ``Lie-Poisson equations".
They are a system of ODE's on $\g^*$.     We will
call the quotient map $T^*G \to (T^*G)/ G = \g^*$ the reduction map.
(Warning: This map is not the reduction  map  of symplectic reduction.)

\subsubsection{Momentum Map}
\label{momentum_map} The left-action of $G$ on itself,
lifted to $T^*G$  has for its momentum map
the map $J: T^* G \to \g^*$ of   {\em right trivialization}.
In terms of our left-trivialized identification
$J(g, M) = Ad_{g^{-1}} ^* M$ where $Ad_g ^* : \g^* \to \g^*$
denotes the dual of the adjoint representation $Ad_g$ of $G$ on $\g$.
The left-invariance of $H$ implies that each  integral curve for the
Hamiltonian vector field $X$, i.e. the geodesics,
when viewed as curves in  the cotangent bundle,  lies within a    constant level set of $J$.

Each individual   constant level-set $J^{-1} (\mu)$ is the image   of a
right-invariant one-form $G \to T^*G$, and as such is a  copy of $G$ in $T^* G$.
 The projection of such a level set onto $\g^*$ by the
 reduction map yields as image the
  co-adjoint orbit through $\mu$ , thus:  $\pi(J^{-1} (\mu)) = G \cdot \mu$
  where $G \cdot \mu = \{ M:  M = Ad_{g} ^* \mu, g \in G \} \subset \g^*$.
  Since the integral curves in $T^*G$ lie on level sets of $J$,
  the   integral curves of the reduced dynamics
  lie on such co-adjoint orbits.

\subsubsection{Unreducing}
\label{unreducing}
Let  $G_{\mu}$ denote the isotropy group of $\mu \in \g^*$
under the co-adjoint action.  As smooth $G$-spaces  we have $\pi(J^{-1} (\mu) ) = G \cdot \mu=  G/ G_{\mu}$ ,
and the projection of $J^{-1} (\mu) \to \pi(J^{-1} (\mu))$
is isomorphic to the canonical bundle projection $G \to G/G_{\mu}$
with fiber $G_{\mu}$.   When $G$ is compact then for
generic $\mu$ we have that
 $G_{\mu}  \cong T$, where $T$ is the maximal torus $T$ of $G$ and the rank $r$ of $G$ is the dimension of $T$.
 If  the typical integral curves  $C$ for the reduced dynamics are
 closed curves $C \subset G \cdot \mu \subset  \g^*$, then the integral curves for  the original dynamics
 sit on manifolds $\pi^{-1} (C) \cap J^{-1} (\mu)$ which is a $T$-bundle over the circle $C$.
 In our particular situation this bundle will be trivial, so that
 $\pi^{-1}(C)\cap J^{-1-}(\mu)$ is itself a torus of one more dimension than $T$.

\subsubsection{Casimirs}
\label{Casimirs}  A Casimir on $\g^*$ is a smooth function
such that  for all smooth functions $h$ on $\g^*$ we have
that $\{C, h \} = 0$.  The values of a Casimir stay constant
on the solutions to the reduced equation.   For $G$  compact  with maximal torus $T$
the algebra of  Casimirs is  functionally generated  by $r$
polynomial generators, these generators being polynomials invariant under
the co-adjoint action.    The common level set  $C_1 = c_1, \ldots,  C_r = c_r$
 of these $r$ Casimirs
is, for generic values of the constants $c_i$, a co-adjoint orbit $G \cdot \mu$ for
  which $G_{\mu} = T$.

 \subsection{The case of Lorentzian metrics on $SU(2) \times S^1$}

 The   Hamiltonian for the  Fefferman metric (eq. \ref{Hamn})
 computed from step 1 is that of  a left-invariant Lorentzian metric on
 $G = SU(2) \times S^1$.   We specialize the discussion of the last few paragraphs
 to this situation.  Then the   dual of the Lie algebra of   $G$ splits as $\g^* = \R^3 \times \R$.
 The   $\R^3$ factor acts like the well-known angular momentum from physics.
    The coordinates $M_1, M_2, M_3, P$ appearing in eq. (\ref{Hamn}) are linear coordinates
 on $\g^* = \R^3 \times \R$. .    Their Lie-Poisson  brackets are
 $$\{ M_1, M_2  \} = - M_3 , \hskip.2cm   \{ M_3, M_1  \} = - M_2, \hskip.2cm  \{ M_2, M_3  \} = - M_1 $$
 together with
 $$\{M_i, P \} = 0.$$
  The rank $r$ of $G$ is $2$.   The algebra of Casimirs is
 generated by 
 $$P \hbox{ and }   K = M_1 ^2 + M_2 ^2 + M_3 ^2 \hskip .3cm \hbox{ (Casimirs) } $$
 Using    momentum scaling, we can split the analysis of the reduced geodesic flow  into two cases,
 $P =0$, and $P =1$.

\subsubsection{{\bf Case 1: $P = 0$}}  We will see that our Hamiltonian equations for this first case are easily solved.
 The reduced dynamics will be trivial:  $M_1 = M_2 = 0$,
 $M_3 = const$.  Up on $G$, the corresponding geodesics
 are left translates of the  one-parameter subgroup corresponding
 to the third direction.\\\\

\subsubsection{\bf Case 2:  $P = 1$} When $P=1$ we
 have for our  Hamiltonian the function $H(M, 1)$ on
 $\R^3 = \R^3 \times \{1 \} \subset \g^*$.   We are only interested
 in the light-like geodesics, which means we will set
 $H(M, 1) = 0$.  This defines a paraboloid in $\R^3$.
 The integral curves for  the reduced dynamics   lie on
 the intersections of this paraboloid with the spheres $K = r_0 ^2$.
 These  intersections   typically consist of one or two closed curves,
 which are the closed  integral curves of the reduced dynamics.

 \subsubsection{Co-adjoint action and identifications}
 \label{coadjoint}
The co-adjoint action of $G$ on $\g^* = \R^3 \times \R$ acts
trivially on the $\R$ factor, since that corresponds to the Abelian factor $S^1$.
 The $\R^3$ factor of $\g^*$  is identified with both $su(2)$ and $su(2)^*$
 and the identification is such that the co-adjoint (or adjoint) action corresponds
 to the standard action of $SO(3)$ on $\R^3$ by way of composition
 with the 2:1 cover $SU(2) \to SO(3)$.   (The $S^1$ factor of $G$
 acts trivially on $\R^3$.)  Under this identification,
 the co-isotropy subgroup $SU(2)_L \subset SU(2)$
 of a  non-zero  vector $L \in \R^3$ consists of the one-parameter subgroup
 generated by $L$, and in $SO(3)$ to rotations about the axis $L$.

\subsubsection{Unreducing}
\label{unreducing2}

 The momentum map $J: T^* G \to \R^3 \times \R$   splits into
 $$ J = (L, J_0) = ((L_1, L_2, L_3), J_0) \hskip .2cm  \hbox{ with } J_0 = P.$$
 The fact that $J_0 = P$ is the $\R$ component of   $J$ is a reflection of the  triviality
 of the  co-adjoint action on the $\R$ factor of  $\g^* =  \R^3 \times \R$.

The solution curves back up on $T^* G$  corresponding  to
a given reduced solution curve $C$ lie on  submanifolds
$J^{-1} (\mu) \cap \pi^{-1} (C)$.  The value of $\mu = (L, P)$
is constrained by the co-adjoint orbit on which $C$ lives.
This constraint is simply $K = \Sigma L_i ^2$.
Only the case $K \ne 0$ is interesting.  Then
the isotropy $G_{\mu}$ is one of  the   maximal torii
$G_{\mu} = SU(2)_L  \times S^1 = S^1 \times S^1 \subset SU(2) \times S^1$.
The first $S^1$ factor  is the  circle $SU(2)_L$   as in the paragraph \ref{coadjoint}.
 It follows
  from the dicussion of (\ref{unreducing}) that
  $J^{-1} (\mu) \cap \pi^{-1} (C)$  is a  $G_{\mu} = S^1 \times S^1$  bundle over $C$.  We
 also saw in (\ref{momentum_map}) that
$J^{-1} (\mu)  \cong G = S^3 \times S^1$.  The projection $\pi$
restricted to  $J^{-1} (\mu)$ is the composition
$S^3 \times S^1 \to S^3 \to S^2 \subset \R^3 \times \{ P = 1 \}$
where the last map is the Hopf fibration.  The Hopf fibration is trivial
over $S^2 \setminus \{ P  \}$ for any point $P \in S^2$.
It follows that $J^{-1} (\mu) \cap \pi^{-1} (C)$   in  isomorphic to a three torus,
$T^3$.  One factor of this
 three-torus is the $S^1$ factor of $SU(2) \times S^1$,
 and corresponds to the extra  angle $\gamma$ we add when constructing
 the circle bundle on which
   Fefferman's  metric lives.   We project out  this angle when
 forming the chains.  Thus the chains lie on two-torii
   $T^2 \subset SU(2)$.   One angle of the two-torus corresponds
 to a coordinate around a curve $C$ in the reduced dynamics.
 The other angle is generated by the circle $SU(2)_L \subset SU(2)$.
 \section{ The reduced Fefferman dynamics.}
\label{phase_portrait}
 \subsection{ The case $P =0$}

When $P= 0$ we see that  $H=\frac{1}{2} (a M_{1}^{2}+\frac{1}{a}
M_{2}^{2})$. Since $H = 0$  we have that $M_1=M_2=0$ along
light-like solutions with $P = 0$. From the constancy of the Casimir
$K$ it follows that $M_3 = const. $ also, so that the reduced
solution is a constant curve.   Generally speaking, for  a
left-invariant metric on a Lie group $G$, the  geodesics in $G$
which  correspond  to  a constant solution  $M(t) = const. = M_*$ of
the reduced equations consist of  the one-parameter subgroup $exp(t
\xi)$ and its left translates $g exp(t \xi)$, where $\bI \xi = M_*$
and $\bI $ is the ``inertial tensor'', i.e. the index lowering
operator corresponding to the metric at the identity.  In our case
$\bI$ maps the $e_3$ axis to the $M_3$ axis, so that the
corresponding geodesic is the  1-parameter subgroup $exp(te_{3})$
and its translations $g exp(t e_3)$.   (More accurately,
 $\bI^{-1} (0,0, M_3, 0)$ is a linear combination of $e_3$
and the basis vector $\frac{\partial}{\partial \gamma}$.  We   project out
the angle $\gamma$  to form the chain corresponding to a light-like geodesic,
so these chains are  indeed generated by $e_3$.)  These $P = 0$  chains are precisely
 circles of    the Hopf fibration   $S^{3} = SU(2)  \to S^2 = SU(2)/S^1$,
 where the $S^1$ is generated by $e_3$ and
 acts by right multiplication.\\\\
\textbf{Note:} Since $-e_3=R$ is the Reeb field for our contact form these chains are the orbits of the Reeb field. It remains to determine whether or not all chains are orbits of Reeb fields.

\subsection{The case  $P = 1$.}

\begin{figure}
\centering
\includegraphics[width=5.5in,height=5.5in]{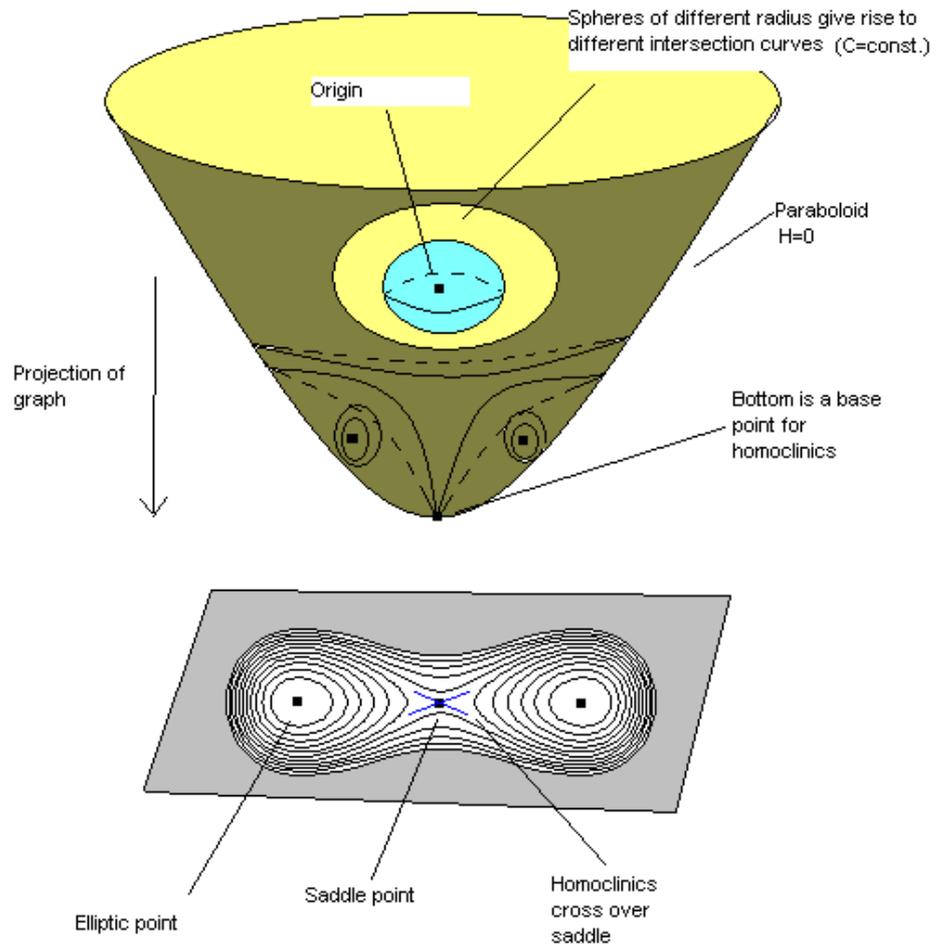}
\caption{The  intersections of  the spheres $K = const.$ with the
paraboloid $H =0$ and their projections to the $M_1 M_2$ plane.}
\label{fig:1}
\end{figure}

\begin{figure}
\centering
\includegraphics[height=2.5in,width=3in]{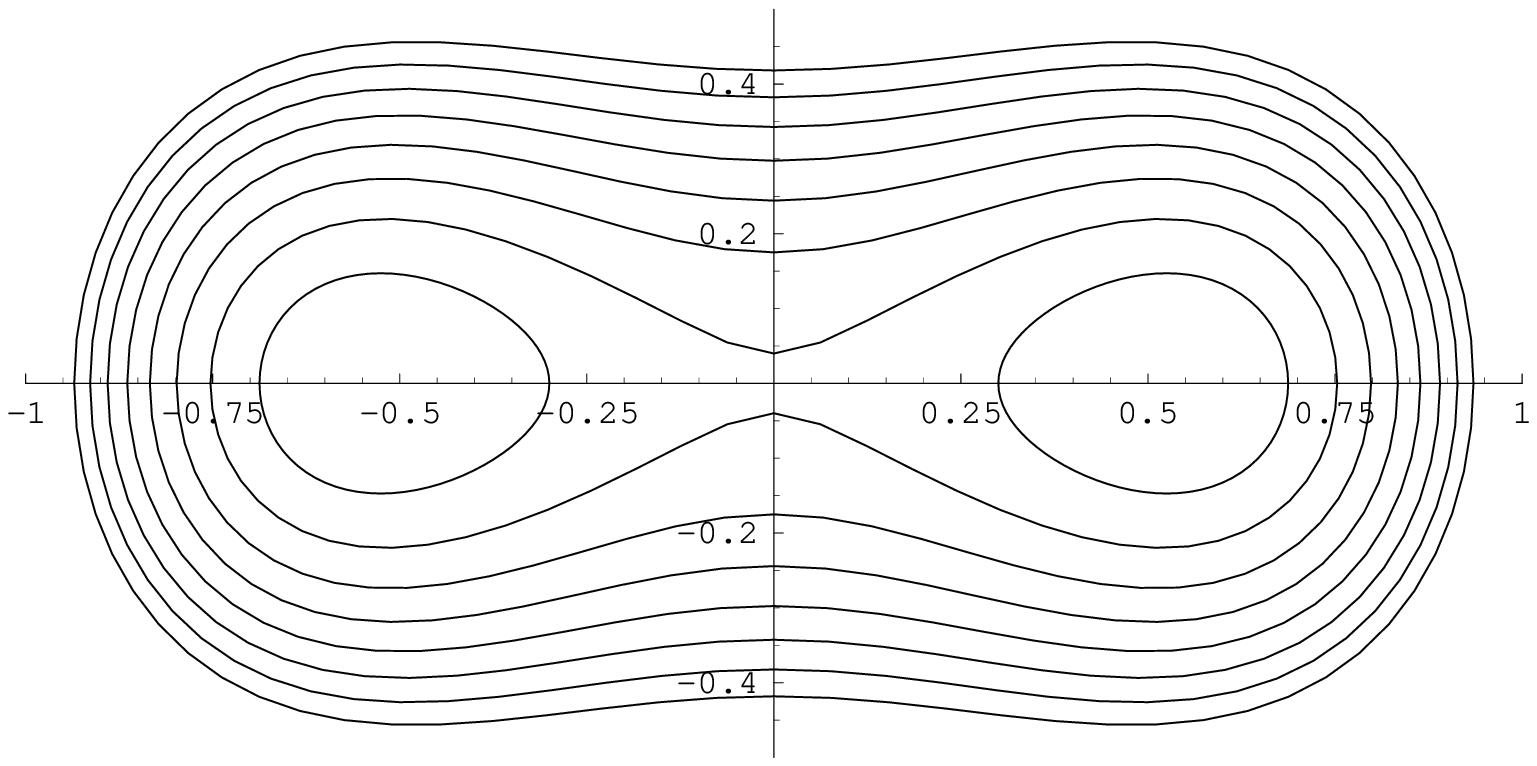}
\caption{Phase portrait for the reduced equations obtained via the projection
in figure 1,   for  $a$ = 2 , $P$ = 1,  and $K$ small.}
\label{fig:2}
\end{figure}

Set $P =1$  in $H$ to   get
$$H_a(M_1,M_2,M_3;1)=\frac{1}{2}(a M_1 ^{2} +\frac{1}{a} M_2 ^{2} - 3M_3 - c(a) )$$
where we have set
$$c(a) =  - \frac{9}{8}(a + \frac{1}{a}).$$
Recall that we are only interested in the solutions for which $H =0$.
The surface $H = 0$ is a {\em paraboloid} which we can express as the graph of a function
of $M_1, M_2$:
\beq
 \{ H = 0 \} = \{ (M_1, M_2, M_3): M_3=\frac{1}{3}(a M_1 ^{2} +\frac{1}{a} M_2 ^{2}  - c(a)) \}
 \label{paraboloid}
 \eeq
The solution curves must also
 lie on   level sets of $K = M_1 ^2 + M_2 ^2 + M_3 ^2$.   In other words,
the solution curves are formed by the
intersection of the paraboloid $H = 0$ with the spheres $K = r_0 ^2$.
See figure \ref{fig:1}.
These intersection curves  are  easily understood by using $M_1, M_2$
as coordinates on the paraboloids, i.e. by projecting the
paraboloid onto the $M_1 M_2$ plane.  They are depicted in figure \ref{fig:2}.

Eq. (\ref{paraboloid}) yields  $M_3$ in terms of $M_1$ and $M_2$ on
the paraboloid. Plug this expression for   $M_3$ into   $K$  to find
that  on the paraboloid
$$K =  (1 -\frac{2}{9} c(a) a ) M_1 ^2 + (1 -\frac{2}{9} \frac{c(a)}{a} ) M_2 ^2  +
 \frac{1}{9} ( a M_1 ^2 + \frac{1}{a} M_2 ^2 )^2 + c(a) ^2.$$
 For $a$ close to $1$ the coefficients of the  quadratic terms, $M_1 ^2$
 and $M_2 ^2$  are   positive, and
 close to $1/2$.   The only critical point for $K$
 is the origin and is a nondegenerate minimum.  It follows from a basic argument in Morse theory
 that all the intersection curves are closed curves, circling  the origin.
 As $a$ increases   the sign of the coefficient  in front of
 the $M_1 ^2$ term  eventually crosses $0$ and becomes negative.
 This happens when $1 -\frac{2}{9} c(a) a  = 0$
 which works out to $a = \sqrt{3}$. After that the origin becomes
 a saddle point for $K$, and the level set of $K$ passing through
 the origin has the shape of a figure 8, with the cross at the origin.
 Inside each lobe of the eight is a new critical point.
  See figure \ref{fig:2} below.
This change as $a$ crosses past $\sqrt{3}$ is an instance of
 what is known as a ``Hamiltonian pitchfork bifurcation'' or ``Hamiltonian
 figure eight'' bifurcation among specialists in Hamiltonian bifurcation theory.

 To re-iterate: for $1 < a < \sqrt{3}$ all reduced solution curves are closed
 and surround the origin.  For $a > \sqrt{3}$ the origin becomes a saddle point,
and the  level set of $K$ passing through the origin
consists of   three
solution curves: the origin itself which is now an   unstable equilibrium,
  and two homoclinic orbits
corresponding to the two lobes of the eight. Being homoclinic to
the unstable equilibrium, it takes an infinite time to traverse either one of
these homoclinic lobes.

 The situation is symmetric as $a$ decreases, with the bifurcation occurring
 at $a = 1/ \sqrt{3}$.  This is as it must be, from the discrete symmetry alluded
 to in Proposition \ref{normalform}, 
  $a \mapsto 1/a$, $M_1 \mapsto M_2, M_2 \mapsto M_1$.

 \section{Step 4:  Berry phase and unreducing.}
 \label{phase}

 As per the discussion in (\ref{unreducing2}),
   associated to each choice of closed
 solution curve
 $C \subset  \R^3 \times \{ 1 \}$
 and each choice  $\mu \ne 0$
 of   momentum,  we have a family of chains
 which lie on  a fixed
 two torus
 $T^2 = T^2 (C; \mu) \subset T^* S^3$.   Our question is :
 are the chains on this   $T^2$   closed?
 The Fefferman dynamics restricted to $T^2$ is
 that of   linear flow on a torus.
   Let $\phi$
 be a choice of angular variable around $C$,
 which we call the base angle.
 Let  $\theta$ be the other  angle of the torus, which we call the `vertical angle'
 chosen so that the projection $T^2 \to C$ is $(\phi, \theta) \mapsto \phi$.
 We take  both angles   defined mod $2 \pi$.
 As we traverse    the chain,   every time that
  the base angle   $\phi$ varies   from $0$ to   $2 \pi$,
 (which is to  say we travel  once around   $C$)     the
 vertical  angle $\theta$ will have varied by some amount
 $\Delta \theta$.  The  amount
 $\Delta \theta$ does not depend on the choice of
 chain within $T^2$.   If $\Delta \theta$ is a rational multiple of $2 \pi$
 then the chains in $T^2$ are all closed.  If $\Delta \theta$ is  an irrational multiple of
 $2 \pi$,
 then none of the chains in $T^2$  close up, and we have the case of quasi-periodic chains
 corresponding to irrational flow on   $T^2$.

Without loss of generality we can suppose that   $\mu = r_0 e_3$ where $e_3$
denotes the final element of  the  standard basis   of  $su(2) ^* = \R^3$.
For why we can assume this without loss of generality refer to subsection
\ref{coadjoint} above. In this case
$K = r_0 ^2$ and this  fixing  of $K$almost fixes the reduced curve $C$.
(See the second paragraph in the proof of the proposition immediately below for
details.)   Remembering the modulus parameter $a$,
we see that
$$\Delta \theta = \Delta \theta (K, a).$$
Since the dynamical system defined by
the Fefferman metric depends analytically on
initial conditions  and on the parameter $a$,
we see that  $\Delta \theta (K, a)$ is an analytic function
of $a$ and $K$.  {\it  It follows that in order to prove theorem 1,
all we need to do is show that for a single value of $a$, the function
$K \mapsto \Delta \theta (K, a)$ is non-constant.}  We see that
in order to prove Theorem 1 it only remains to prove:

\begin{proposition}
\label{phaseVarn}
For $a > \sqrt{3}$ the  function
$K \mapsto \Delta \theta (K, a)$ is non-constant.
\end{proposition}

{\bf Proof of Proposition.}

Fix  $a > \sqrt{3}$.    Consider the
value $K = c(a)^2$ corresponding to the homoclinic figure eight
through the origin in the   $M_1 M_2$ plane.
We will show that
\beq
\lim_{K \to c(a)^2 _- } \Delta \theta (K, a) = +\infty.
\label{limit}
\eeq
and that for $K$ slightly less than $c(a)^2$ the value of
$\Delta \theta (K, a)$ is finite.  It follows that
the function $K \mapsto \Delta \theta (K, a)$ varies, as required.

Let $m(a)$ denote the absolute minimum of $K$ on the
paraboloid.  The minimum is achieved
at two points, the  elliptic fixed points inside each    lobe of the homoclinic eight.
For values of $r_0 ^2$ between $m(a)$ and $c(a)^2$ the
level set $K = r_0 ^2$ consists of two disjoint closed curves $C_1, C_2$, one inside each
lobe of the eight.
These two curves  are related  by the reflection
$(M_1, M_2) \mapsto (-M_1, M_2)$.  The entire
dynamics is invariant under this reflection,   so that  the
value of $\Delta \theta$ on $C_1$
equals   its value on $C_2$.     (The two components
are traversed in the same sense.)  Consequently   $\Delta \theta (K, a)$ is well-defined
and finite
for $m(a) < K < c(a)^2$,  being equal to the common  value of $\Delta \theta (C_i)$.

In what follows we arbitrarily fix one of the
two components of $K = r_0 ^2$ and call it   $C$.

The key to establishing the limit (\ref{limit}) is a Berry phase
formula for $\Delta \theta$ which   mimics earlier work of one of us
(\cite{rmontC}).    The formula expresses $\Delta \theta$
as the sum of two integrals:
\beq
\Delta \theta (K, a)  = \hbox{dynamic} + \hbox{ geometric }
\label{BerryPhase}
\eeq
where
$$\hbox{ dynamic } ={1 \over {\sqrt{K}}}  \int_0 ^T f dt  $$
and
$$\hbox{ geometric} = - (\hbox{oriented solid angle}).$$
Both the dynamic and the geometric terms can be expressed as line integrals around
  $C$.
In the dynamic term, $T = T(K)$ is the period of the curve $C$,
and where
\beq
f = \frac{1}{2} [a M_1 (t)  ^2 + \frac{1}{a} M_2 (t) ^2 + c(a) ].
\label{integrand}
\eeq
The integral is done around the     projection of the curve $C$
to the $M_1 M_2$.
  The time    $t$ is the time
parameter occuring in  the reduced equations, which is the same as the
geodesic time.
In the second formula, the oriented solid angle is the standard oriented  solid angle
enclosed by a closed  curve such as $C$ in space.  The
 absolute value of an oriented   solid angle is always  bounded by
$4 \pi$.  On the other hand,
$ \frac{1}{\sqrt{K}} f >   \frac{1}{2 \sqrt{K}} c(a)$. Consequently,   if we let  the curve $C$ approach the
lobe of the homoclinic
orbit  which contains it,   then its  period $T(K) $ tends to $\infty$.
We now see that   the   dynamic term
of     eq. (\ref{BerryPhase})
tends to $+\infty$.
 Thus, the  corollary is proved  once we have  established the validity of the
Berry phase type formula (\ref{BerryPhase}).

\subsection{Proof of Berry phase formula}

 We  begin the proof of eq. (\ref{BerryPhase}) by  recalling
  and summarizing our situation,
  and applying    the discussion of    (\ref{unreducing2})
  for relating the
  reduced dynamics to  dynamics in   $T^*(SU(2) \times S^1)$ and curves in  $T^* SU(2)$.
 We have  fixed  $J = (L, P)$ to equal  the value $\mu = (r_0 e_3, 1) \in  R^3 \times \R$
 where  $r_0 \ne 0$.    The values of the
 Casimirs  which characterize our reduced curve $C$ are then  $K = r_0 ^2$, and $P =1$.
 The Fefferman light-like geodesics $C_F$ associated
 to $C$ and our choice of $\mu$ must   lie on the manifold
 $J^{-1} (\mu) \cap \pi ^{-1} (C)$ which
   is a three-torus inside $T^*(SU(2) \times S^1)$.
 Project  this three torus  into   $T^* S^3$ via the product structure induced projection:
 $pr_2: T^* (S^3 \times S^1) = T^* S^3 \times T^*S^1 \to T^* S^3$ and in this way
 arrive at a  two-torus  $X(C) = pr_2 (J^{-1} (\mu)) \cap \pi ^{-1} (C)  \subset T^* SU(2) \times \{1 \}$ which
 projects onto $C$ via the canonical projection $T^* (SU(2)) \times \{1 \} \to \R^3 \times \{1 \}$.
 We will  soon need that
 $X(C) \subset  L^{-1} (r_0 e_3) \times \{1 \}$ which follows from  the fact that $J = (L, P)$
 so that $pr_2 (J^{-1} (\mu)) = L^{-1} (r_0 e_3) \times \{1 \}$.
 The canonical projection just refered to is that of the   quotient map  $T^* (SU(2)) \to \R^3$
 for  the (lifted) left action of $SU(2)$ on itself.   The momentum map
 associated to this map is $L$.
   We will also use that the  canonical projection, $T^* (SU(2)) \to \R^3$,
 restricted to level sets of $L$, corresponds to  symplectic reduction for $T^*SU(2)$.
 The chains  $ch$ associated to the reduced solution $C$ and our
 choice of momentum axis $e_3$ lie in the
 two-torus $X(C)$.  To coordinatize $X(C)$ choose any global section
 $\hat C: C \to  X(C)$ and let $\phi$ be an angular coordinate
 around $C$ so that   $\hat C$ is a closed curve in  $X(C)$ parameterized by $\phi$
 and projecting onto $C$.
 Now act on $\hat C$  by the one-parameter subgroup
 $exp(\theta e_3) = SU(2)_{L}$.   Then any point of
 $X(C)$ can be written as $exp(\theta e_3) \cdot \hat C (\phi) \in X(C)$
 where  $\theta, \phi$ are global angular coordinates.
 (The multiplication ``$ \cdot $''
 of  ``$exp(\theta e_3) \cdot \hat C (\phi)$''  denotes the action of the group element
 $exp(\theta e_3)  \in SU(2)$
 on $T^*SU(2)$ by cotangent lift.)

 Every cotangent bundle $T^* Q$ is endowed with a canonical one-form.
 Let $\Theta$ be the canonical one-form on $T^* SU(2)$.
 Our Berry phase formula (\ref{BerryPhase}) will be  proved by applying Stoke's theorem
 to the integral of $\Theta$ around a well-chosen closed curve $c$ in
 $X(C)$.

 This curve $c \subset X(C) \subset T^*SU(2) \times \{1 \}$ is the concatenation of
 two curves.  One curve is  any one of the chains $ch$ corresponding to $C$ -- which is to say,
 the projection  by $pr_2$ of any one of the Fefferman geodesics $C_F  \subset J^{-1} (\mu) \cap \pi ^{-1} (C)$.  We parameterize $ch$ by the Fefferman dynamical time, $0 \le  t \le T$
 making sure to stop when, upon projection,we have gone once round $C$,
 so that $C(0) = C(T)$.
 Having gone once round $C$, we must have $ch(T) = exp(\Delta \theta e_3) \cdot c(0)$.
 The {\em holonomy} $\Delta \theta$ is the angle we are trying to compute.
 For the other curve $c_{group}$ we simply move backwards in the
 group direction to close up the curve:
 $c_{group}(s) =   exp(-s e_3) \cdot ch(T)$.
 Our   curve $c$ is then the concatenation  $+$ of these
 two smooth curves:
 $$c = c_{group}  +  ch.$$

 The curve $c$ is a closed curve lying in the two-torus  $X(C)$.
 Not all closed curves in the two-torus bound discs,
 but $X(C) \subset  L^{-1} (r_0 e_3) \times \{1 \}  \cong SU(2)$
 which is simply connected, so that $c$ does bound a disc
   $ \tilde D \subset L^{-1} (r_0 e_3) \times \{1 \}$.
 Apply Stoke's formula:   \beq
 \int_{\tilde D} d \Theta  = \int_{c_{group}}  \Theta + \int_{ch} \Theta.
 \label{Stokes}
 \eeq
 The proof of (\ref{BerryPhase}) proceeds by evaluating  each term in eq (\ref{Stokes}) separately.

 Write $S^2$ for the two-sphere $K = r_0^2$, $P =1$.  Write
   $\pi_{r_0}  :  L^{-1} (r_0 e_3)  \to  S^2$ for the restriction
 of the canonical reduction map  $\pi: T^* SU(2)  \times \{ P = 1 \} \to  \R^3 \times \{1\}$.
  Under $\pi_{r_0}$ the disc $\tilde D$ projects onto
 a  topological disc $D \subset S^2$  which bounds our reduced  curve $C$.  $S^2$ is the symplectic reduced space of
 $T^*SU(2)$
 by the left action of $SU(2)$, reduced  at the value $L = r_0 e_3$.  A  basic result from symplectic reduction, essentially its
 definition, asserts that  as a  symplectic reduced space  $S^2$
 is endowed with a 2-form $\omega_{r_0}$ (the reduced symplectic form)
 defined by $\pi_{r_0} ^* \omega_{r_0} = i^* (- d \Theta)$,
 where $i:   L^{-1} (r_0 e_3)  \to T^* SU(2)$ is the inclusion.
 Let $d \Omega$ denote  the unique rotationally invariant
 two-form on the two sphere, normalized so that its integral over
 the entire sphere is $4 \pi$.   (The form $d \Omega$  is not closed, but the notation is
standard, and  suggestively helpful, so we use it.)
 It is well-known that $\omega_{r_0} = -r_0 d \Omega$,
 which is to say, that
 $$r_0 ( \pi_{r_0} ^* ( d \Omega) )= i_{r_0} ^* (d \Theta).$$
 (See \cite{Foundations} for the standard ``high-tech'' computation, and
 \cite{rmontC} for an elementary computation of this well-known fact.)
 Thus
 \beq
 \int _{\tilde D} ( d \Theta )  =   \int_D r_0 d \Omega = r_0 (\hbox{ solid angle enclosed by  } C )
 \label{area}
 \eeq
 It is worth noting that this area is a signed area, positive or negative
 depending on the orientation of the bounding curve $C$ of $D$.

  It follows from the definition of the momentum map
 on the cotangent bundle  that $\Theta ( {d \over {ds}} (exp(s e_3) (p)) = r_0$
 for any point $p \in L^{-1} (r_0 e_3)$.
 It follows that
 $$\Theta = r_0 d \theta  \hbox{ along }  c_{group},$$
 and thus
 \beq
 \int_{c_{group}}  \Theta = - r_0 \Delta \theta.
 \label{group}
 \eeq
 where
 the minus sign arises because in travelling along $c_{group}$
 we moved backwards in the $e_3$-direction.

 It remains to compute $\int_{ch} \Theta$.
For this computation we will have to work on $T^* (SU(2) \times S^1 ) $.
There we have the canonical one form
 \beq
 \Theta_F = \Theta + P d \gamma.
 \label{decomposition}
 \eeq
 Now relative to   any coordinates $x^a$ for $SU(2) \times S^1$,
 where $p_a$ are the corresponding momentum coordinates   we have
 $$\Theta_F = \Sigma p_a dx^a. $$
 Plugging in along one  of the light-like Fefferman geodesics
 and using the metric relation $p_a = \Sigma g_{ab} \dot x ^a $
 where $g_{ab}$ are the metric components
 we see that
 $$\Theta_F ( \dot C_F (t)) = 2 H = 0$$
 where the last equality arises because  the Fefferman geodesic is light-like.
 Since
 $pr_2 \circ C_F = ch$
 where   $pr_2: T^*SU(2) \times T^* S^1 \to T^* SU(2)$ is the projection,
 we have, from (\ref{decomposition}),
 $$\Theta ( {d \over {dt}}  ch  ) = - P \dot \gamma = - \dot \gamma,$$
where we used $P =1$.  It follows that
 $$\int_{ch} \Theta = - \int_0 ^T  \dot \gamma  dt $$.
 Now $\dot \gamma = \frac{\partial H} {\partial P}$.
 Referring back to the equation for the Hamiltonian, and remembering
 that we set $P =1$ after differentiating we see that
  $$\dot \gamma = -\frac{3}{2} M_3 - c(a).$$
  Now using the formula for $M_3$ in terms of $M_1, M_2$ and a bit of algebra
  we see that
  $$- \dot \gamma =  f, $$
where $f$ is as in the eq. (\ref{integrand}).  Thus: \beq \int_{ch}
\Theta = \int_0 ^T f dt. \label{dynamic} \eeq

 Putting together the pieces (\ref{area}), (\ref{group}), (\ref{dynamic})
  into Stokes' formula (\ref{Stokes}) and some algebra yields
  the Berry phase formula (\ref{BerryPhase}). QED

  \vskip 2cm
  \appendix{{\bf APPENDICES}}

 \section{ The dynamics when $a =1$.}
 \label{app1}

The chains for the standard structure on $S^3$  are       formed by  intersecting
$S^3 \subset \C^2$ with complex lines in $\C^2$.  See \cite{Goldman}.  In this appendix we
verify that the Fefferman metric description of chains when $a =1$
yields these circles.

The key to our verification   is   the observation that when $a=1$
the Fefferman  Hamiltonian (26) splits into two commuting pieces $H
= H_0 - H_1$ with $\{H_0, H_1 \} = 0$. This observation and the
following method of computation is the same one which led to
explicit formulae for \sR geodesic flows in chapter 11 of
\cite{montA}, formulae identical to that of Lemma 1 below. We have
$H_0 = {1\over 2} K = {1 \over 2} (M_1 ^2 + M_2 ^2 + M_3 ^2)$ and
$H_1 = \frac{1}{2}(M_3 - \frac{3}{2} P) ^2$.   Since the two
Hamiltonians commute, their flows up on the cotangent bundles
commute. This observation leads to the explicit formula for the
chains through the identity: \beq ch(t) = exp\left[ t(M_1e_1 + M_2
e_2 + M_3e_3)\right] exp\left[- t(M_3-\frac{3}{2}P)e_3)\right]
\label{chains} \eeq The $M_i, P$ are constants which  satisfy the $H
= 0$ condition
$$(M_1^{2}+M_2^{2}+M_{3}^{2})=(M_{3} -\frac{3}{2} P)^{2}.$$
In this formula (\ref{chains}) for the chains, the first factor corresponds to
the flow of $H_0$, whose integral curves correspond to one-parameter
subgroups in $SU(2)$, and the second factor corresponds to the projection
to $SU(2)$ of solutions to the Hamilton's equation for $-H_1$.

To verify that the chains computed via Fefferman's metric are the circles
described above
we  use  two lemmas from linear algebra.\\

\textbf{Lemma 1.} ( circles in SU(2)) \textit{Every geometric  circle in $SU(2) = S^3$  through the identity can be parameterized as  $\gamma(t)=exp(\alpha t)exp(-\beta t)$ where $\alpha,\beta \in su(2)$ are
Lie algebra elements of the same length. } \\

\textbf{Lemma 2.} \textit{When $\beta = c e_3$ as in eq. (\ref{chains}) then these circles sit on complex lines.}\\

{\bf Remark.}  The condition $|\alpha| = |\beta|$ in lemma 1 is a $1:1$ resonance
condition.

\vskip .2cm

The proofs rely on identifying  the quaternions $\mathbb{H}$ with
$\mathbb{C}^{2}$ and hence  the group of unit quaternions  with
$SU(2)$ and  $S^3$.  Since the contact plane is annihilated by
$\omega_3$, and is to correspond with the $T_x S^3 \cap
\mathbb{J}T_x S^3$, we must take the  identification $\C^2 \cong
\bH^2$ such that the complex structure on $\C^2$    corresponds to
right multiplication by $k$, where  $k$ is to  correspond to $e_3$
in $su(2)$.

\textbf{Proof of   lemma 1}.  In a Euclidean vector space, (such as
$\bH$) the  circles are described by $c(t) = P + r (\cos(\omega t)
e_1 + \sin(\omega t) e_2)$ where $P$ is the center of the circle,
$r$ its radius,  and  where $e_1, e_2$ are an \on basis for the
plane through $P$ containing the circle. Now use the fact that for a
unit quaternion $n$ we have $exp(n t) = \cos(t) 1 + \sin (t) n$.
Thus $\gamma(t)$ of lemma 1 is equal to $(\cos(t)  + \sin(t)
\alpha)(\cos(t)  -\sin(t) \beta)$. Algebra and trigonometry
identities   yield
$$\gamma (t) = \frac{1}{2} [ (1 - \alpha \beta) + \cos(2 t) (1 + \alpha \beta) + \sin(2t)( \alpha - \beta)]$$
which we can rewrite as
$$\gamma (t) = P + \cos(2 t) v + \sin(2t) w, $$
with $P =\frac{1}{2}(1 - \alpha \beta)$, $v =  \frac{1}{2}  (1 +
\alpha \beta)$  and $w =  \frac{1}{2} ( \alpha - \beta)$. It remains
to show that $v$ and $w$ have the same length and are orthogonal.
Using $\bar \alpha = - \alpha$ and remembering that $\alpha$ is unit
length we see that we have $v = - \alpha w$ and so indeed $|v| =
|w|$.  Their   common length is the radius $r$ of the circle. Since
the Euclidean inner product is given by $Re(v \bar w)$ the fact that
$v = - \alpha w$ also shows that $v$ and $w$ are orthogonal. QED

\textbf{Proof of   lemma 2. }   Let $v, w$ be as in the proof of lemma 1.
 We must show that
the real 2-plane spanned by $v$ and $w$ is a complex line
when $\beta = k$.   Recall that under our identification of
$\C^2$ with $\bH$  the complex structure corresponds
to multiplication on the right by $k$.  Now compute   $wk = v$,
to see that the span of $v$ and $w$ is  indeed a complex line. QED

 \section{ Relation to the Rossi example. }
 \label{app2}

Rossi \cite{rossi} constructed a much-cited example
 of a family of non-embeddable CR-structures on $S^3$.
 The purpose of this appendix is to show that Rossi's family
 is isomorphic to  our   left-invariant CR family with $a \ne 1$. This isomorphism
 is well-known to experts. We include it here for completeness.     We use the description of CR manifolds
 to be found in the remark towards the beginning of section \ref{circle bundle}.
 In that construction a CR structure is defined as the span of complex
 vector field.  Let $Z$ be the complex vector field corresponding to  the
 standard CR structure.  In terms of our left invariant frame,
 $Z = e_1 - i e_2$.
  Then Rossi's perturbed
 CR structure is defined by
 $$Z_{\mu} = Z - \mu \bar Z$$
 with $\mu$ a real parameter.
 On the other hand, we saw (again,  eq. \ref{holomorphic}) that
 our left-invariant CR structures correspond to the span
 of
 $$Z_a = e_1 - \frac{i}{a}e_2. $$
 Set  $a=1+\epsilon$ and expand out:
$Z_a  = e_1 - i(1+\epsilon)e_2=e_1 - i e_2 - i\epsilon e_{2}=Z +\frac{1}{2}\epsilon(Z-\bar{Z})$.
Upon rescaling $Z_a$ by dividing by  $(1 + \frac{1}{2}\epsilon)$ we see that
$span(Z_a) = Span( Z - \mu(\epsilon) \bar{Z}$,
where $\mu(\epsilon)=\frac{\frac{1}{2}\epsilon}{1+\frac{1}{2}\epsilon}$.
This shows that the left-invariant structure for $a$ corresponds to
Rossi's structure for $\mu = \mu(\epsilon)$.  \\\\

The important facts concerning Rossi's structures for $\mu \ne 0$
is that every CR-function for one of these structures on
$S^3$ is \textit{even} with respect to the antipodal map $(x,y,z) \mapsto(-x,-y,-z)$.
We recommend Burns' \cite{Burns} for the proof.
This forced evenness  s implies that  there is no CR embedding
of  our left-invariant structures for $a \ne 1$ into    $\C^n$ for any $n$.
The structures  do however, have  explicit $2:1$ immersions into $\C^3$
which  can  be found in Rossi.  See also  Burns (\cite{Burns}) or  Falbel \cite{Falbel}.
Upon taking the quotient by the antipodal map each  $a \ne 1$
structure induces a left-invariant  CR structure on
$\R P ^3 = SO(3)$ which does  embed into $\C^3$.
This embedded image bounds a domain within  an explicit Stein manifold
  $S \subset \C^3$.

  \vskip 1cm

 {\sc Open Problem. }  [Dan Burns] Find a synthetic construction of the
    chains for the left-invariant structures, in the spirit of the
construction
    of the chains for the standard structure, but using a family
of complex
   curves in $S$ in place of the straight lines used to construct
the chains for the standard structure.

\section{Acknowledgements.}   We would
like to thank John Lee for explaining the connection between the
Rossi example and the left-invariant structures as detailed in
Appendix 2, and Dan Burns for e-mail conversations for further
considerations concerning Appendix 2, and for the open problem, and
for listening critically to an early version of our results. We
would like to acknowledge encouragement and helpful conversations
from  Gil Bor of CIMAT, Jie Qing (UCSC), and Robin Graham (U. of
Washington). We would especially like to thank Gil Bor for crucial
help regarding the reduced dynamics, and CIMAT (Guanajuato, Mexico)
for support during this time. The junior author also thanks AGEP for
some financial aid in many instances of the project.

\end{document}